\def\balpha{\mbox{\boldmath $\alpha$}}

\def\bmu{\mbox{\boldmath $\mu$}}

\def\bsigma{\mbox{\boldmath $\sigma$}}

\def\Bsigma{{\bf \Sigma}}

\def\Ba{{\bf A}}
\def\Bb{{\bf B}}

\def\Bi{{\bf I}}

\def\Bs{{\bf S}}

\def\Bu{{\bf U}}
\def\Bv{{\bf V}}

\def\Bx{{\bf X}}

\def\bolda{{\bf a}}
\def\boldb{{\bf b}}
\def\boldc{{\bf c}}

\documentclass[10pt]{article}

\newcommand{\myeq}[1]{(\ref{#1})}

\begin{document}

\title{ Asymptotic distribution for the proportional covariance model}
\author{Myung Geun Kim\footnote{ E-mail: mgkim@seowon.ac.kr}}
\date{}

\maketitle

\begin{abstract}
Asymptotic distribution  for the proportional covariance model under multivariate normal distributions is derived. To
this end, the parametrization of the common covariance matrix by its Cholesky root is adopted. The derivations are made
in three steps. First, the asymptotic distribution of the maximum likelihood estimators of the proportionality
coefficients and the Cholesky inverse root of the common covariance matrix is derived  by finding the information
matrix and its inverse. Next,  the asymptotic distributions for the case of the Cholesky root of the common covariance
matrix  and finally for the case of the common covariance matrix itself are derived using the multivariate
$\delta$-method. As an application of the asymptotic distribution derived here, a hypothesis for homogeneity of
covariance matrices is considered.
\end{abstract}

 \vspace{.2in}


\section{Introduction}

\hspace{\parindent}
Let $\Bx_{k}$ be a generic $p$-variate random vector having  a multivariate
normal distribution
$N_{p}(\bmu_{k},\Bsigma_{k})$ for each $k=1,...,K,$ where the $\bmu_{k}$ are
column vectors of the $p$-dimensional
Euclidean space and the $\Bsigma_{k}$ are $p$ by $p$ positive definite covariance
matrices.
 Let the $\Bs_{k}$
 be the mutually independent unbiased estimators of the $\Bsigma_{k}$,
each for a sample of size $N_{k}$ and let $n_{k}=N_{k}-1$.
 Then the $n_{k}\Bs_{k}$ are
independently distributed according to  Wishart distribution
$W_{p}(\Bsigma_{k},n_{k})$.

The hypothesis that $K$ covariance matrices are proportional to each other
can be expressed as
\begin{equation} \label{eq:1d1}
 H_{p}:\Bsigma_{k}=c_{k}\Bsigma_{1}
\mbox{~~~~for all  k = 2, ... ,K}
\end{equation}
 where the $c_{k}$ are
 unknown positive proportionality coefficients.
In the two sample case, estimation for proportional covariance models was studied by
Khatri(1967),  Guttman, Kim and Olkin (1985) and Rao (1983).  The proportional
covariance model was adopted in classification by 0wen (1984). An algorithm to find
the maximum likelihood estimators was given independently by 0wen (1984), Manly and
Rainer (1987) and Eriksen (1987).  The convergence of the algorithm and the
uniqueness of maximum likelihood estimates were proved by Eriksen (1987). Using a
parametrization of the spectral decomposition of the common covariance matrix, Flury
(1986) considered maximum likelihood estimation of the proportionality coefficients
and the common covariance matrix, and suggested an iterative method for performing
computations. Under the same parametrization as in Flury (1986),  Boente et al.
(2007, 2009) studied influence diagnostics and robust estimation for proportional
covariance models.

The Cholesky decomposition theorem states that there exists a lower
triangular matrix $\Ba$ with positive diagonal elements such that
$$\Bsigma_{1} = \Ba\Ba^{T},$$
and this expression is unique. We use $T$ to denote a transpose of a matrix. Then
the model \myeq{eq:1d1} is parametrized as
\begin{equation} \label{eq:1d2}
\Bsigma_{k} = c_{k}\Ba\Ba^{T}~~~(1 \leq k \leq K).
\end{equation}
 For  notational convenience, it
is assumed that $c_{1}$ is identical with 1. The parametrization \myeq{eq:1d2} can
be considered as the special case of the model treated by J\"{o}reskog (1971) and
the easily accessible LISREL (linear structural relationships) program is available
for getting the maximum likelihood estimates of the model parameters (J\"{o}reskog
and S\"{o}rbom, 1996).

Let $\Bb^{T} = \Ba^{-1}$. The matrix $\Bb$ is then an upper triangular matrix with
positive diagonal elements. Then the parametrization \myeq{eq:1d2} can be written as
\begin{equation} \label{eq:1d3}
\Bsigma_{k}^{-1} = \frac{1}{c_{k}}\Bb\Bb^{T}~~~(1 \leq k \leq K).
\end{equation}

In this work, we derive the asymptotic distribution of the maximum likelihood
estimators of the proportionality coefficients and the common covariance matrix,
using the well-known properties of the maximum likelihood estimators under
regularity conditions (Rao, 1973, Chapter 6). The procedure for doing so is as
follows. Firstly, the information matrix for the parameters $c_{k}$ and $\Bb$ in
\myeq{eq:1d3} is derived. Secondly, the asymptotic covariance matrix for the maximum
likelihood estimators of the parameters $c_{k}$ and $\Bb$ in \myeq{eq:1d2} is found
by inverting the corresponding information matrix. Thirdly, we derive the asymptotic
covariance matrix for the maximum likelihood estimators of the parameters $c_{k}$
and $\Ba$ in \myeq{eq:1d2}, using the multivariate $\delta$-method. Finally, the
asymptotic covariance matrix for the maximum likelihood estimators of the
proportionality coefficients and the common covariance matrix is derived using the
multivariate $\delta$-method. Finally,  a hypothesis for homogeneity of covariance
matrices is considered as an application of the asymptotic distribution derived
here.

\section{The information matrix for the parameters $c_{k}$ and $\Bb$}

\hspace{\parindent} In this section, we derive the information matrix for the
$c_{k}$ and $\Bb$. The likelihood function of the $\Bsigma_{k}$ given the $\Bs_{k}$
is
$$
 L(\Bsigma_{1},...,\Bsigma_{K}) = C \times \prod_{k=1}^{K}
[ \mid \Bsigma_{k}
\mid^{-\frac{n_{k}}{2}} \mbox{exp}\{-\frac{n_{k}}{2}\mbox{tr}
(\Bsigma_{k}^{-1}\Bs_{k})\}],
$$
where $C$ is a constant not depending on the parameters. Since we are interested in
covariance matrices and a covariance matrix is location invariant, there is no
restriction on mean vectors.

Under the parametrization \myeq{eq:1d3}, we have
\begin{eqnarray*}
 \mid \Bsigma_{k}^{-1} \mid & = & c_{k}^{-p} \prod_{i=1}^{p}b_{ii}^{2} \\
 \mbox{tr}(\Bsigma_{k}^{-1}\Bs_{k}) & = &
 \frac{1}{c_{k}}\sum_{i=1}^{p}\boldb_{i}^{T}\Bs_{k}\boldb_{i},
\end{eqnarray*}
where $b_{ij}$ denotes the $(i,j)$th element of $\Bb$ and $\boldb_{i}$ the
$i$th column of $\Bb$. Thus the log-likelihood function, ignoring
constant term, becomes
$$ l(\boldc,\Bb) = n_{+}\sum_{i=1}^{p}\mbox{log}(b_{ii})
- \sum_{k=1}^{K}\frac{n_{k}}{2}\{p\mbox{log}(c_{k}) +
 \frac{1}{c_{k}}\sum_{i=1}^{p}\boldb_{i}^{T}\Bs_{k}\boldb_{i}\},$$
 where $n_{+} = \sum_{k=1}^{K}n_{k}$ and $\boldc = (c_{2},...,c_{K})^{T}$.
 It is easily shown that all the expectations of
first derivatives of $l(\boldc,\Bb)$ with respect to the parameters are zero,
 and  the expectations of  second derivatives of the log-likelihood function
 are as follows:
 \begin{eqnarray*}
 - E\{\frac{\partial^{2}l(\boldc,\Bb)}{\partial c_{k}^{2}}\} & = &
 \frac{pn_{k}}{2c_{k}^{2}}~~~~(2 \leq k \leq K)\\
 - E\{\frac{\partial^{2}l(\boldc,\Bb)}{\partial b_{ii}^{2}}\} & = &
n_{+}(a_{ii}^{2} + \bolda_{(i)}^{T}\bolda_{(i)})~~~~(1 \leq i \leq p)\\
 - E\{\frac{\partial^{2}l(\boldc,\Bb)}{\partial b_{ii} \partial c_{k}}\} & = &
 - \frac{n_{k}a_{ii}}{c_{k}}~~~~(2 \leq k \leq K,~1 \leq i \leq p)\\
 - E\{\frac{\partial^{2}l(\boldc,\Bb)}{\partial b_{ji} \partial c_{k}}\} & = & 0
 ~~~~(2 \leq k \leq K,~j < i)\\
 - E\{\frac{\partial^{2}l(\boldc,\Bb)}{\partial b_{ji} \partial b_{j'i}}\} & = &
 n_{+}\bolda_{(j)}^{T}\bolda_{(j')}~~~~(1 \leq j, j' < i),
 \end{eqnarray*}
where $a_{ij}$ denotes the $(i,j)$th element of $\Ba$ and $\bolda_{(i)}$ the
column vector formed by the elements in the $i$th row of $\Ba$,
and all the other cases have zero expectations.

Let $\Bi (\boldc,\Bb)$ be the  information matrix  for the parameters $\boldc$ and
$\Bb$, and let it be partitioned as
$$ \Bi (\boldc,\Bb) = \pmatrix{
\Bi_{11} (\boldc,\Bb) & \Bi_{12} (\boldc,\Bb)\cr
\Bi_{21} (\boldc,\Bb) & \Bi_{22} (\boldc,\Bb)\cr},$$
such that $\Bi_{11} (\boldc,\Bb)$ is the information matrix for $\boldc$ and
  $\Bi_{22} (\boldc,\Bb)$ for the nonzero elements of $\Bb$.
Let $r_{k} = n_{k}/n_{+}.$ Then $\Bi_{11} (\boldc,\Bb)$ is a
 diagonal matrix of order $K-1$:
$$
\frac{p}{2}\mbox{diag}(\frac{r_{2}}{c_{2}^{2}},...,\frac{r_{K}}{c_{K}^{2}}).
$$

We are interested in nonzero elements of $\Bb$, $b_{ij}~(i < j)$ and
the number of nonzero elements of $\Bb$ is $p(p+1)/2$.
Let $\boldb_{i.j}$~$= (b_{1i},...,b_{ji})^{T}.$ Then we can notice that
$\boldb_{i.1} = b_{1i}$ and $\boldb_{i.p} = \boldb_{i}.$
Let $\balpha$~$= (r_{2}/c_{2},...,r_{K}/c_{K})^{T}.$ We write as ${\bf 1}_{p.i}$
the unit column vector
of dimension $p$ having one in the $i$th position and
zeroes in the other positions.
Note that ${\bf 1}_{1.1} = 1.$
Then the $(K-1)$ by $p(p+1)/2$ matrix $\Bi_{12} (\boldc,\Bb)$ is
$$
\bordermatrix{
     & \boldb_{1.1}^{T}  & \boldb_{2.2}^{T} & \ldots & \boldb_{p.p}^{T} \cr
\boldc & - a_{11}\balpha & - a_{22}\balpha {\bf 1}_{2.2}^{T}
& \ldots & - a_{pp}\balpha {\bf 1}_{p.p}^{T} }.
$$

Let $\Ba_{i}$ be the leading principal submatrix of $\Ba$ having order $i$.
Note that $\Ba_{1} = a_{11}$ and $\Ba_{p} = \Ba.$
Then the $p(p+1)/2$ by $p(p+1)/2$ matrix  $\Bi_{22} (\boldc,\Bb)$ is
 a block diagonal matrix whose $i$th block element for $\boldb_{i.i}$
is given by
$$
a_{ii}^{2} {\bf 1}_{i.i} {\bf 1}_{i.i}^{T} + \Ba_{i}\Ba_{i}^{T}.
$$
This $i$th block element is of order $i$.

\section{Asymptotic distribution of the maximum likelihood estimators
of the $c_{k}$ and $\Bb$}

\hspace{\parindent}
The asymptotic distribution of the maximum likelihood estimators of
$\boldc$ and $\Bb$ is characterized by its asymptotic covariance matrix,
and the asymptotic covariance matrix is just the inverse of the
corresponding information matrix.

Let $\Bv (\boldc,\Bb)$ be the asymptotic covariance matrix for the maximum
likelihood estimators of the parameters $\boldc$ and the nonzero elements
of $\Bb$, and let it be partitioned as
$$ \Bv (\boldc,\Bb) = \pmatrix{
\Bv_{11} (\boldc,\Bb) & \Bv_{12} (\boldc,\Bb)\cr
\Bv_{21} (\boldc,\Bb) & \Bv_{22} (\boldc,\Bb)\cr},$$
such that $\Bv_{11} (\boldc,\Bb)$ is the covariance matrix
for the maximum likelihood estimator of  $\boldc$ and
  $\Bv_{22} (\boldc,\Bb)$ for the maximum likelihood estimators of
  the nonzero elements of $\Bb$.

\vspace{.15in}
The following lemma gives the inverse of $\Bi_{22} (\boldc,\Bb)$ and
enables us to get $\Bv_{11} (\boldc,\Bb)$,

{\bf Lemma 3.1.}~~~For each $i = 1,...,p,$ the inverse of
$a_{ii}^{2} {\bf 1}_{i.i} {\bf 1}_{i.i}^{T} +$~$\Ba_{i}\Ba_{i}^{T}$ is
$$\Bb_{i}\Bb_{i}^{T} - \frac{1}{2}\boldb_{i.i}\boldb_{i.i}^{T}$$
where $\Bb_{i}$ is the leading principal submatrix of $\Bb$ having
dimension $i$.

\vspace{.15in}

\noindent Lemma 3.1 is easily proved using the result (Sherman-Morrison
formula) by Bartlett (1951) and
the fact that $\Bb_{i}^{T} = \Ba_{i}^{-1}.$

\vspace{.15in}

{\bf Lemma 3.2.}~~~The matrix $\Bi_{11} (\boldc,\Bb)~-$~
$\Bi_{12}(\boldc,\Bb)\Bi_{22}(\boldc,\Bb)^{-1}\Bi_{21} (\boldc,\Bb)$,
$\Bu_{11}(\boldc,\Bb)$ say, of dimension $K-1$ is given by
$$
\frac{p}{2}\{\mbox{diag}(\frac{r_{2}}{c_{2}^{2}},...,\frac{r_{K}}{c_{K}^{2}})
- \balpha\balpha^{T}\}.
$$

{\bf Proof of Lemma 3.2.}~~~The $i$th block element of
$\Bi_{12}(\boldc,\Bb)\Bi_{22}(\boldc,\Bb)^{-1}$ is a $(K-1)$ by $i$ matrix
for the pair $\boldc$ and $\boldb_{i.i}$ and is given by
$$ - \frac{1}{2}\balpha\boldb_{i.i}^{T}.$$
Then $\Bi_{12}(\boldc,\Bb)\Bi_{22}(\boldc,\Bb)^{-1}\Bi_{21} (\boldc,\Bb)$
is given by
$$\frac{p}{2}\balpha\balpha^{T}$$
and therefore the proof is completed.

\vspace{.15in}

The inverse of $\Bu_{11}(\boldc,\Bb)$ is just $\Bv_{11} (\boldc,\Bb)$
and is stated in the following lemma.

\vspace{.15in}

{\bf Lemma 3.3.}~~~The $(K-1)$ by $(K-1)$ matrix $\Bv_{11} (\boldc,\Bb)$
is given by
$$
\frac{2}{p}\{\mbox{diag}(\frac{c_{2}^{2}}{r_{2}},...,\frac{c_{K}^{2}}{r_{K}})
+ \frac{1}{r_{1}}\boldc\boldc^{T}\}.
$$

\vspace{.15in}

\noindent The matrix $\Bv_{11} (\boldc,\Bb)$ is easily derived using the result
by Bartlett (1951).

\vspace{.15in}

The matrix $\Bi_{12}(\boldc,\Bb)\Bi_{22}(\boldc,\Bb)^{-1}$ is derived in the
proof of Lemma 3.2 and
 $\Bv_{12} (\boldc,\Bb)~=$ ~ $- \Bv_{11} (\boldc,\Bb)
\Bi_{12} (\boldc,\Bb)\Bi_{22} (\boldc,\Bb)^{-1}$.
Hence we get easily the following lemma.

{\bf Lemma 3.4.}~~~The $(K-1)$ by $p(p+1)/2$ matrix $\Bv_{12} (\boldc,\Bb)$ is
given by the following matrix multiplied by $(1/pr_{1})$:
$$
\bordermatrix{
     & \boldb_{1.1}^{T} & \boldb_{2.2}^{T} & \ldots & \boldb_{p.p}^{T} \cr
\boldc & \boldc \boldb_{1.1}^{T} & \boldc \boldb_{2.2}^{T}
& \ldots & \boldc \boldb_{p.p}^{T} }.
$$

\vspace{.15in}

Since $\Bv_{22} (\boldc,\Bb)~=$~$\Bi_{22} (\boldc,\Bb)^{-1} -$~
$\Bi_{22} (\boldc,\Bb)^{-1}$$\Bi_{21} (\boldc,\Bb)$$\Bv_{12} (\boldc,\Bb)$,
 Lemma 3.1, the proof of Lemma 3.2 and Lemma 3.4 yield $\Bv_{22} (\boldc,\Bb)$.

{\bf Lemma 3.5.}~~~
The $p(p+1)/2$ by $p(p+1)/2$ matrix  $\Bv_{22} (\boldc,\Bb)$ is given by
$$
\bordermatrix{
     & \boldb_{1.1}^{T} & \boldb_{2.2}^{T} & \ldots & \boldb_{p.p}^{T} \cr
\boldb_{1.1} & \Bb_{1}\Bb_{1}^{T} + d\boldb_{1.1}\boldb_{1.1}^{T} &
e\boldb_{1.1}\boldb_{2.2}^{T} & \ldots & e\boldb_{1.1}\boldb_{p.p}^{T} \cr
\boldb_{2.2} & e\boldb_{2.2}\boldb_{1.1}^{T} &
\Bb_{2}\Bb_{2}^{T} + d\boldb_{2.2}\boldb_{2.2}^{T} &
\ldots & e\boldb_{2.2}\boldb_{p.p}^{T} \cr
\vdots & \vdots & \vdots & \ddots & \vdots \cr
\boldb_{p.p} & e\boldb_{p.p}\boldb_{1.1}^{T} &
e\boldb_{p.p}\boldb_{2.2}^{T} & \ldots &
\Bb_{p}\Bb_{p}^{T} + d\boldb_{p.p}\boldb_{p.p}^{T}},
$$
where $d = \{1-(1+p)r_{1}\}/2pr_{1}$ and $e = (1-r_{1})/2pr_{1}$.

\vspace{.15in} We use hereafter the {\em hat} notation to denote the maximum
likelihood estimator.

\vspace{.15in}
 The general properties of the maximum likelihood estimators
 yield the asymptotic distribution of $\hat{\boldc}$ and $\hat{\Bb}$
 under the proportional covariance model.

\vspace{.15in}

{\bf Theorem 3.1.}~~~As $n_{+}$ tends to infinity with each $r_{k}$
held fixed, the asymptotic distribution of
$$ \sqrt{n_{+}}(\hat{\boldc} - \boldc, \hat{\Bb} - \Bb)$$
is a multivariate normal with mean zero and covariance matrix
$\Bv(\boldc,\Bb)$.

\vspace{.15in}

Lemma 3.3 shows that the asymptotic distribution of
$\sqrt{n_{+}}(\hat{c}_{2} - c_{2},...,\hat{c}_{K} - c_{K})$
is free of $\Bsigma_{1}$. Guttman et al. (1985) derived the asymptotic
distribution
of $\sqrt{n_{+}}(\hat{c}_{2} - c_{2})$ when $K=2$ and Lemma 3.3 reduces
to Theorem 4.1 of Guttman et al. (1985).

When $r_{1} = 1$, that is, the single-population case, the coefficients
$d$ and $e$ become $- 1/2$ and zero, respectively. Hence Lemma 3.5
shows that the $\hat{\boldb}_{i.i}$~$(i=1,...,p)$ are asymptotically
independent.
By the invariance property of the maximum likelihood estimator,
the Cholesky inverse root of $\Bs_{1}$ is just $\hat{\Bb}$,
that is,
$\Bs_{1}^{-1} = \hat{\Bb}\hat{\Bb}^{T}.$ Thus Lemma 3.5 gives the asymptotic
distribution of the Cholesky inverse root of the sample covariance matrix.

\section{Asymptotic distribution of the maximum likelihood
estimators of  the $c_{k}$ and $\Ba$}

\hspace{\parindent}
The asymptotic distribution of the maximum likelihood estimators
of the $c_{k}$ and $\Ba$ is derived from that of the $c_{k}$ and
$\Bb$, using the multivariate $\delta$-method.

Let $\Bv (\boldc,\Ba)$ be the asymptotic covariance matrix for the maximum
likelihood estimators of the parameters $\boldc$ and $\Ba$. We write
$\Bv (\boldc,\Ba)$ as a partitioned form:
$$ \Bv (\boldc,\Ba) = \pmatrix{
\Bv_{11} (\boldc,\Ba) & \Bv_{12} (\boldc,\Ba)\cr
\Bv_{21} (\boldc,\Ba) & \Bv_{22} (\boldc,\Ba)\cr},$$
such that $\Bv_{11} (\boldc,\Ba)$ is the covariance matrix
for the maximum likelihood estimator of  $\boldc$ and
  $\Bv_{22} (\boldc,\Ba)$ for the maximum likelihood estimators of
  the nonzero elements of $\Ba$.

Since $\Bb^{T} = \Ba^{-1}$, we have the following lemma.

\vspace{.15in}

{\bf Lemma 4.1.}~~~
The partial derivative of $a_{ji}$ with respect to $b_{hg}$ is
$$
\frac{\partial  a_{ji}}{\partial b_{hg}} =
\left\{ \begin{array}{cl}
- a_{jg}a_{hi} & (i \leq h \leq g \leq j)\\
0  &  \mbox{otherwise.}
\end{array} \right.
$$

\vspace{.15in}

Let $\bolda_{i.j} = (a_{1i},...,a_{ji})^{T}$ and
$\bolda_{i(-j)}$~$= (a_{j+1,i},...,a_{pi})^{T}$. Note that $\bolda_{i(-0)}$ and
$\bolda_{i.p}$ are equivalent to $\bolda_{i}$, the $i$th column of $\Ba$
and that $\bolda_{i.1} = a_{1i}$.
Then the Jacobian, $J(\Ba,\Bb)$, of $\Ba$ with respect to $\Bb$ is
$$
\bordermatrix{
&\boldb_{1.1}^{T}&\boldb_{2.2}^{T}&\boldb_{3.3}^{T}&\ldots&\boldb_{p.p}^{T}\cr
\bolda_{1} & -\bolda_{1}\bolda_{1.1}^{T} & -\bolda_{2}\bolda_{1.2}^{T}
& -\bolda_{3}\bolda_{1.3}^{T} & \ldots & -\bolda_{p}\bolda_{1.p}^{T} \cr
\bolda_{2(-1)} & {\bf 0}  & -\bolda_{2(-1)}\bolda_{2.2}^{T}
& -\bolda_{3(-1)}\bolda_{2.3}^{T} & \ldots
& -\bolda_{p(-1)}\bolda_{2.p}^{T} \cr
\bolda_{3(-2)} & {\bf 0}  & {\bf 0} &
-\bolda_{3(-2)}\bolda_{3.3}^{T} &  \ldots & -\bolda_{p(-2)}\bolda_{3.p}^{T} \cr
\vdots & \vdots & \vdots & \vdots & \ddots & \vdots \cr
\bolda_{p(-p+1)} & {\bf 0}  & {\bf 0} &
{\bf 0} & \ldots &
-\bolda_{p(-p+1)}\bolda_{p.p}^{T} }.
$$
Hence the Jacobian $J(\Ba,\Bb)$ is a block upper triangular matrix.

\vspace{.15in}
Under the transformation $\Bb^{T} = \Ba^{-1}$, the asymptotic distribution
of  $\hat{\boldc}$ is not changed and in fact $\Bv_{11} (\boldc,\Ba) =$~
$\Bv_{11} (\boldc,\Bb)$.

Since $\Bv_{12} (\boldc,\Ba)$~$= \Bv_{12} (\boldc,\Bb)J(\Ba,\Bb)^{T}$,
we get the following lemma.

\vspace{.15in}

{\bf Lemma 4.2.}~~~The $(K-1)$ by $p(p+1)/2$ matrix $\Bv_{12} (\boldc,\Ba)$ is
given by the following matrix multiplied by $-1/pr_{1}$:
$$
\bordermatrix{
     & \bolda_{1}^{T}& \bolda_{2(-1)}^{T} & \ldots & \bolda_{p(-p+1)}^{T} \cr
\boldc & \boldc \bolda_{1}^{T} & \boldc \bolda_{2(-1)}^{T}
& \ldots & \boldc \bolda_{p(-p+1)}^{T} }.
$$

{\bf Proof of Lemma 4.2.}~~~By Lemma 3.4, $\Bv_{12} (\boldc,\Ba)$ has as
its $i$th block element for the covariance matrix of $\hat{\boldc}$ and
$\hat{\bolda}_{i(-i+1)}$
$$ - \frac{\boldc}{pr_{1}}
\sum_{j=i}^{p}\boldb_{j.j}^{T}\bolda_{i.j}\bolda_{j(-j+1)}^{T}
~=~ - \frac{1}{pr_{1}}\boldc\bolda_{i(-i+1)}^{T}$$
because $\Bb^{T}\Ba = \Bi_{p}$.
\vspace{.15in}

Let $\Ba_{-i}$ be the matrix formed by deleting the first $i$ columns and the first
$i$ rows of $\Ba$ simultaneously. Note that $\Ba_{-0} = \Ba$ and $\Ba_{-p+1} =
\bolda_{p(-p+1)} = a_{pp}$.

\vspace{.15in}

{\bf Lemma 4.3.}~~~
The $p(p+1)/2$ by $p(p+1)/2$ matrix  $\Bv_{22} (\boldc,\Ba)$ is given by
$$
\bordermatrix{
     & \bolda_{1}^{T}&\bolda_{2(-1)}^{T}& \ldots &\bolda_{p(-p+1)}^{T} \cr
\bolda_{1} & \Ba\Ba^{T} + d\bolda_{1}\bolda_{1}^{T} &
e\bolda_{1}\bolda_{2(-1)}^{T} & \ldots & e\bolda_{1}\bolda_{p(-p+1)}^{T} \cr
\bolda_{2(-1)} & e\bolda_{2(-1)}\bolda_{1}^{T} &
\Ba_{-1}\Ba_{-1}^{T} + d\bolda_{2(-1)}\bolda_{2(-1)}^{T} &
\ldots & e\bolda_{2(-1)}\bolda_{p(-p+1)}^{T} \cr
\vdots & \vdots & \vdots & \ddots & \vdots \cr
\bolda_{p(-p+1)} & e\bolda_{p(-p+1)}\bolda_{1}^{T} &
e\bolda_{p(-p+1)}\bolda_{2(-1)}^{T} & \ldots &
(1+d)\bolda_{p(-p+1)}\bolda_{p(-p+1)}^{T}},
$$
where $d$ and $e$ are  defined as in Lemma 3.5.

{\bf Proof of Lemma 4.3.}~~~The $(i,j)$th block element of $\Bv_{22} (\boldc,\Ba)$
is the covariance matrix of $\hat{\bolda}_{i(-i+1)}$ and
$\hat{\bolda}_{j(-j+1)}$, and is given by
$$
\sum_{l=i}^{p}\sum_{m=j}^{p}
\bolda_{l(-i+1)}\bolda_{i.l}^{T}\Bv_{22}[l,m]
\bolda_{j.m}\bolda_{m(-j+1)}^{T}~~~~(i \leq j),
$$
where $\Bv_{22}[l,m]$ denotes the $(l,m)$th block element of
$\Bv_{22} (\boldc,\Bb)$ which is the covariance matrix of $\hat{\boldb}_{l.l}$
and $\hat{\boldb}_{m.m}$. Let $\delta_{ij}$ be the Kronecker delta.
Since $\bolda_{i.l}^{T}\boldb_{l.l}~=$~$\delta_{il}$ and
$\bolda_{i.l}^{T}\Bb_{l}~=$~${\bf 1}_{l.i}^{T}$, the $(i,j)$th block element
of $\Bv_{22} (\boldc,\Ba)$ becomes
$$
e(1-\delta_{ij})\bolda_{i(-i+1)}\bolda_{j(-j+1)}^{T} +
\delta_{ij}(\Ba_{-i+1}\Ba_{-i+1}^{T} + d\bolda_{i(-i+1)}\bolda_{i(-i+1)}^{T}),
$$
which completes the proof.

\vspace{.15in}

Lemmas 3.3, 4.2 and 4.3 yield the asymptotic distribution of $\hat{\boldc}$
and $\hat{\Ba}$ under the proportional covariance model.

\vspace{.15in}

{\bf Theorem 4.1.}~~~The asymptotic distribution of
$$ \sqrt{n_{+}}(\hat{\boldc} - \boldc, \hat{\Ba} - \Ba)$$
is a multivariate normal with mean zero and covariance matrix $\Bv (\boldc,\Ba)$,
as $n_{+}$ tends to infinity with each $r_{k}$ kept fixed.

\vspace{.15in}

Consider the single-population case, $r_{1} = 1$. Since the coefficients $d$
and $e$ become $- 1/2$ and zero, respectively, Lemma 4.3 shows that the
$\hat{\bolda}_{i(-i+1)}$
~$(1 \leq i \leq p)$ are asymptotically independent. In this single-population
case, $\Bs_{1} = \hat{\Ba}\hat{\Ba}^{T}$ and therefore
Lemma 4.3 gives the asymptotic distribution of
 the Cholesky root of the sample covariance matrix
when the underlying distribution is a multivariate one.
The asymptotic distribution of the sample generalized variance being equal to
the product of the $\hat{a}_{ii}^{2}$ is easily
derived from Lemma 4.3 and it can be found also in Theorem 7.5.4 of
Anderson (1984).

\section{Asymptotic distribution of the maximum likelihood estimators
of the $c_{k}$ and $\Bsigma_{1}$}

\hspace{\parindent}
Now, we are in a position to derive the asymptotic distribution of the
maximum likelihood estimators of the proportionality coefficients
and the common covariance matrix, and it will be done using the multivariate
$\delta$-method.

Let $\sigma_{ij}$ be the $(i,j)$th element of $\Bsigma_{1}$.
The common covariance matrix $\Bsigma_{1}$ has $p(p+1)/2$ independent
parameters.
Let $\Bi_{p.i}$ be the $(p-i)$ by $p$ matrix in which the first $i$
columns are zero vectors and the remaining $p-i$ columns form the unit
matrix of order $p-i$. When $i = 0$, $\Bi_{p.0}$ becomes the unit matrix of
order $p$ and is denoted by $\Bi_{p}$.
For the $i$th column $\bsigma_{i}$ of $\Bsigma_{1}$, $\bsigma_{i(-j)}$
is defined similarly to $\bolda_{i(-j)}$.

Since $\sigma_{ij}$~$= \bolda_{(i)}^{T}\bolda_{(j)}$, we get the Jacobian,
$J(\Bsigma_{1},\Ba)$, of $\Bsigma_{1}$ with respect to $\Ba$ as follows:
$$
\bordermatrix{
 & \bolda_{1}^{T}&\bolda_{2(-1)}^{T}&\bolda_{3(-2)}^{T}&\ldots&\bolda_{p(-p+1)}^{T} \cr
 & a_{11}\Bi_{p} &  &  &  &  \cr
\bsigma_{1} & +   &{\bf 0} & {\bf 0} & \ldots & {\bf 0} \cr
 &  \bolda_{1} {\bf 1}_{p.1}^{T} & & & & \cr
 & a_{21}\Bi_{p.1} & a_{22}\Bi_{p-1}  &  &  &  \cr
\bsigma_{2(-1)} &  +  & + & {\bf 0} & \ldots & {\bf 0} \cr
 &  \bolda_{1(-1)} {\bf 1}_{p.2}^{T} &\bolda_{2(-1)} {\bf 1}_{p-1.1}^{T} & & & \cr
 & a_{31}\Bi_{p.2} & a_{32}\Bi_{p-1.1}  & a_{33}\Bi_{p-2} &  &  \cr
\bsigma_{3(-2)} &  +  & + & + & \ldots & {\bf 0} \cr
 &  \bolda_{1(-2)} {\bf 1}_{p.3}^{T} &\bolda_{2(-2)} {\bf 1}_{p-1.2}^{T}
 &\bolda_{3(-2)} {\bf 1}_{p-2.1}^{T} & & \cr
 \vdots & \vdots & \vdots & \vdots & \ddots & \vdots \cr
 & a_{p1}\Bi_{p.p-1} &a_{p2}\Bi_{p-1.p-2}  & a_{p3}\Bi_{p-2.p-3} &
 &a_{pp}  \cr
\bsigma_{p(-p+1)} &  +  & + & + & \ldots & + \cr
 &  \bolda_{1(-p+1)} {\bf 1}_{p.p}^{T} &\bolda_{2(-p+1)} {\bf 1}_{p-1.p-1}^{T}
 &\bolda_{3(-p+1)} {\bf 1}_{p-2.p-2}^{T} &
 &\bolda_{p(-p+1)} }.
$$
Hence the Jacobian $J(\Bsigma_{1},\Ba)$ is a block lower triangular matrix.

\vspace{.15in}

Let $\Bv (\boldc,\Bsigma_{1})$ be the asymptotic covariance matrix
for the maximum
likelihood estimators of the parameters $\boldc$ and $\Bsigma_{1}$, and let it
be partitioned as
$$ \Bv (\boldc,\Bsigma_{1}) = \pmatrix{
\Bv_{11} (\boldc,\Bsigma_{1}) & \Bv_{12} (\boldc,\Bsigma_{1})\cr
\Bv_{21} (\boldc,\Bsigma_{1}) & \Bv_{22} (\boldc,\Bsigma_{1})\cr},$$
such that $\Bv_{11} (\boldc,\Bsigma_{1})$ is the covariance matrix
for the maximum likelihood estimators of $\boldc$ and
  $\Bv_{22} (\boldc,\Bsigma_{1})$ for the maximum likelihood estimator of
the independent parameters of  $\Bsigma_{1}$.

The covariance matrix of $\hat{\boldc}$, $\Bv_{11} (\boldc,\Bsigma_{1})$,
is equivalent to $\Bv_{11} (\boldc,\Bb)$ given in Lemma 3.3.

\vspace{.15in}

We need some useful formulas to get $\Bv_{12} (\boldc,\Bsigma_{1})$ which are stated
in the following lemma which is easily proved.

\vspace{.15in}

{\bf Lemma 5.1.}~~~For $j \leq i$, we have\\
(i)~~~$\Bi_{p-j+1.i-j}\bolda_{j(-j+1)} = \bolda_{j(-i+1)}$\\
(ii)~~$\bolda_{j(-j+1)}^{T}{\bf 1}_{p-j+1.i-j+1} = a_{ij}$\\
(iii)~$\bsigma_{i(-i+1)} = \sum_{j=1}^{i}a_{ij}\bolda_{j(-i+1)}.$\\

\vspace{.15in}

The covariance matrix $\Bv_{12} (\boldc,\Bsigma_{1})$ of $\hat{\boldc}$ and
$\hat{\Bsigma}_{1}$ is given by $\Bv_{12} (\boldc,\Ba)J(\Bsigma_{1},\Ba)^{T}$ as in
the following lemma.

\vspace{.15in}

{\bf Lemma 5.2.}~~~The $(K-1)$ by $p(p+1)/2$ matrix
$\Bv_{12} (\boldc,\Bsigma_{1})$ is
given by the following matrix multiplied by $-2/pr_{1}$:
$$
\bordermatrix{
     & \bsigma_{1}^{T}&\bsigma_{2(-1)}^{T}& \ldots &\bsigma_{p(-p+1)}^{T} \cr
\boldc & \boldc \bsigma_{1}^{T} & \boldc \bsigma_{2(-1)}^{T}
& \ldots & \boldc \bsigma_{p(-p+1)}^{T} }.
$$

{\bf Proof of Lemma 5.2.}~~~The covariance matrix of $\hat{\boldc}$ and
$\hat{\bsigma}_{i(-i+1)}$ which is the $i$th block element of
$\Bv_{12} (\boldc,\Bsigma_{1})$ is given by
$$
- \frac{1}{pr_{1}}\boldc \sum_{j=1}^{i}
\bolda_{j(-j+1)}^{T}
(a_{ij}\Bi_{p-j+1.i-j}^{T} + {\bf 1}_{p-j+1.i-j+1}\bolda_{j(-i+1)}^{T})
$$
which becomes by Lemma 5.1
$$ - \frac{2}{pr_{1}}\boldc \bsigma_{i(-i+1)}^{T}.$$

\vspace{.15in}

Let $\Bsigma_{(i,j)} = (\bsigma_{j(-i+1)},...,\bsigma_{p(-i+1)})$.
Then $\Bsigma_{(i,j)}$ is a $(p-i+1)$ by $(p-j+1)$ matrix
consisting of the elements in the lower right corner of $\Bsigma_{1}$
starting at $\sigma_{ij}$. The $(p-i+1)$ by $(p-j+1)$ matrix $\Ba_{(i,j)}$
is defined similarly to $\Bsigma_{(i,j)}$.

\vspace{.15in}

{\bf Lemma 5.3.}~~~
The $p(p+1)/2$ by $p(p+1)/2$ matrix  $\Bv_{22} (\boldc,\Bsigma_{1})$ has
as its $(i,j)$ block element
$$ \sigma_{ij}\Bsigma_{(i,j)} + 4e\bsigma_{i(-i+1)}\bsigma_{j(-j+1)}^{T} +
\bsigma_{j(-i+1)}\bsigma_{i(-j+1)}^{T}~~~~(i \leq j)$$
which is the covariance matrix of $\hat{\bsigma}_{i(-i+1)}$ and
$\hat{\bsigma}_{j(-j+1)}$,
where $d$ and $e$ are  defined as in Lemma 3.5.

\vspace{.15in}

Since $3+4d = 4e+1$, the covariance matrix of  $\hat{\bsigma}_{i(-i+1)}$
becomes
$$ \sigma_{ii}\Bsigma_{(i,i)} + (3+4d)\bsigma_{i(-i+1)}\bsigma_{i(-i+1)}^{T}.$$

More formulas are needed to prove Lemma 5.3. Let $\bolda_{(i)(-j)}~=$~
$(a_{i,j+1},$...$,a_{ip})^{T}.$ Then we get the following lemma together
with Lemma 5.1 useful for proving Lemma 5.3.

\vspace{.15in}

{\bf Lemma 5.4.}~~~
For $i \leq j$,~$1 \leq m \leq i$ and $1 \leq l \leq j$, we have\\
\vspace{.1in}

\noindent (i)~~~~~${\bf 1}_{p-m+1.i-m+1}^{T}\Ba_{-m+1} = \bolda_{(i)(-m+1)}^{T}$\\
(ii)~~~~$\Bi_{p-m+1.i-m}\Ba_{-m+1} = \Ba_{(i,m)}$\\
(iii)~~~$\bolda_{(i)(-m+1)}^{T}\Ba_{(j,m)}^{T} = \sum_{u=m}^{i}a_{iu}
\bolda_{u(-j+1)}^{T}$\\
(iv)~~~$\bolda_{(i)(-m+1)}^{T}\bolda_{(j)(-m+1)} = \sum_{u=m}^{i}a_{iu}
a_{ju}$\\
(v)~~~~$\bsigma_{i(-i+1)}\bsigma_{j(-j+1)}^{T} =
\sum_{l=1}^{j}\sum_{m=1}^{i}a_{im}a_{jl}\bolda_{m(-i+1)}\bolda_{l(-j+1)}^{T}$\\
\begin{tabbing}
(vi)~~~$\bsigma_{j(-i+1)}\bsigma_{i(-j+1)}^{T}$ \= $=
\sum_{l=1}^{j}\sum_{m=1}^{i}a_{im}a_{jl}\bolda_{l(-i+1)}\bolda_{m(-j+1)}^{T}$\\
\> $= \sum_{m=1}^{i}a_{im}a_{jm}\bolda_{m(-i+1)}\bolda_{m(-j+1)}^{T}$\\
\> $+ \sum_{m=1}^{i}\{a_{jm}\bolda_{m(-i+1)}(\sum_{l=m+1}^{i}a_{il}\bolda_{l(-j+1)}^{T})\}$\\
\> $+ \sum_{m=1}^{i}\sum_{l=m+1}^{j}a_{im}a_{jl}\bolda_{l(-i+1)}\bolda_{m(-j+1)}^{T}$\\
\end{tabbing}
(vii)~~$\Bsigma_{(i,j)} = \sum_{u=1}^{m-1}\bolda_{u(-i+1)}\bolda_{u(-j+1)}^{T}
+ \Ba_{(i,m)}\Ba_{(j,m)}^{T}$\\
\begin{tabbing}
(viii)~$\sigma_{ij}\Bsigma_{(i,j)}$ \= $=
(\sum_{m=1}^{i}a_{im}a_{jm})(\sum_{l=1}^{p}\bolda_{l(-i+1)}\bolda_{l(-j+1)}^{T})$\\
\> $= \sum_{m=1}^{i}a_{im}a_{jm}\Ba_{(i,m)}\Ba_{(j,m)}^{T}$\\
\> $+ \sum_{l=1}^{i-1}(\sum_{m=l+1}^{i}a_{im}a_{jm})\bolda_{l(-i+1)}\bolda_{l(-j+1)}^{T}.$\\
\end{tabbing}

\vspace{.1in}

{\bf Proof of Lemma 5.3.}~~~For $i \leq j$, the $(i,j)$th block element of
$\Bv_{22}(\boldc,\Bsigma_{1})$ which is the covariance matrix of
$\hat{\bsigma}_{i(-i+1)}$ and $\hat{\bsigma}_{j(-j+1)}$ is given by
$$\sum_{l=1}^{j}[\sum_{m=1}^{i}\{(a_{im}\Bi_{p-m+1.i-m} +
\bolda_{m(-i+1)}{\bf 1}_{p-m+1.i-m+1}^{T})\Bv_{22}[m,l]\}$$
\begin{equation} \label{eq:5d1}
(a_{jl}\Bi_{p-l+1.j-l}^{T} + {\bf 1}_{p-l+1.j-l+1}\bolda_{l(-j+1)}^{T})],
\end{equation}
where $\Bv_{22}[m,l]$ denotes the $(m,l)$th block element of $\Bv_{22}(\boldc,\Ba)$.
By Lemmas 5.1 and 5.4, summing the summands in \myeq{eq:5d1} over $1 \leq l = m \leq
i$ gives \begin{eqnarray} \label{eq:5d2}
 &\sum_{m=1}^{i}\{ a_{im}a_{jm}\Ba_{(i,m)}\Ba_{(j,m)}^{T} +
4da_{im}a_{jm}\bolda_{m(-i+1)}\bolda_{m(-j+1)}^{T} \nonumber \\
+ &a_{jm}\bolda_{m(-i+1)}\bolda_{(i)(-m+1)}^{T}\Ba_{(j,m)}^{T} +
a_{im}\Ba_{(i,m)}\bolda_{(j)(-m+1)}\bolda_{m(-j+1)}^{T}\nonumber \\
+ & (\sum_{u=m}^{i}a_{iu}a_{ju})\bolda_{m(-i+1)}\bolda_{m(-j+1)}^{T}\}
\end{eqnarray}
 and taking a summation of the summands in \myeq{eq:5d1} over $ l \neq m$
yields
\begin{equation} \label{eq:5d3}
4e\bsigma_{i(-i+1)}\bsigma_{j(-j+1)}^{T} -
4e\sum_{m=1}^{i}a_{im}a_{jm}\bolda_{m(-i+1)}\bolda_{m(-j+1)}^{T}.
\end{equation}
Applying Lemma 5.4 again to \myeq{eq:5d2} and \myeq{eq:5d3} completes the proof
since $4d + 3 = 1 + 4e$.

\vspace{.15in}

From Lemmas 3.3, 5.2 and 5.3, we get the asymptotic distribution of
 $\hat{\boldc}$ and
$\hat{\Bsigma}_{1}$ under the proportional covariance model \myeq{eq:1d1}.

\vspace{.15in}

{\bf Theorem 5.1.}~~~The statistic
$$ \sqrt{n_{+}}(\hat{\boldc} - \boldc, \hat{\Bsigma}_{1} - \Bsigma_{1})$$
is asymptotically distributed according to a multivariate normal
distribution with mean zero and covariance matrix $\Bv(\boldc,\Bsigma_{1})$
as $n_{+}$ tends to infinity with each $r_{k}$ held fixed.

\vspace{.15in}

The asymptotic distribution of $\hat{\boldc}$ is free of $\Bsigma_{1}$
as can be seen in Lemma 3.3, and that of $\hat{\Bsigma}_{1}$ is
independent of $\boldc$. However, $\hat{\boldc}$ and $\hat{\Bsigma}_{1}$
are not asymptotically independent.

When $r_{1} = 1$, the $(i,j)$th block element of $\Bv_{22}(\boldc,\Bsigma_{1})$
becomes
\begin{equation} \label{eq:5d4}
\sigma_{ij}\Bsigma_{(i,j)} + \bsigma_{j(-i+1)}\bsigma_{i(-j+1)}^{T}
\end{equation}
because $e$ is zero. In this case, the maximum likelihood estimator of $\Bsigma_{1}$
becomes $\Bs_{1}$, derived from the Wishart distribution, which is an unbiased
estimator of $\Bsigma_{1}$. However, the maximum likelihood estimator of
$\Bsigma_{1}$ based on the multivariate normal distribution is in fact
$(n_{1}/N_{1})\Bs_{1}$. Since $\Bs_{1}$ and $(n_{1}/N_{1})\Bs_{1}$ have the same
limiting distribution, this case $(r_{1} = 1)$ yields the asymptotic distribution of
the maximum likelihood estimator of $\Bsigma_{1}$. We can easily check if
\myeq{eq:5d4} is identical with (15) in Anderson (1984, p.82).

\section{A hypothesis for homogeneity of covariance matrices}

\hspace{\parindent}
As an application of Theorem 5.1, consider a hypothesis for  homogeneity
of covariance matrices
\begin{equation} \label{eq:6d1}
H_{0}: c_{2} = . . . = c_{K} =1.
\end{equation}
The inverse $\Bu_{11} (\boldc,\Bb)$ of $\Bv_{11} (\boldc,\Bsigma_{1})$ is needed to
test the hypothesis \myeq{eq:6d1}, which is given in Lemma 3.2.

The asymptotic distribution of
\begin{displaymath}
n_{+}(\hat{c}_{2} - 1,...,\hat{c}_{K} - 1)\Bu_{11} (\boldc,\Bb)(\hat{c}_{2} - 1,...,
\hat{c}_{K} - 1)^{T}
\end{displaymath}
under \myeq{eq:6d1} is a chi-squared distribution with degrees of freedom $K-1$.
However, $\Bu_{11} (\boldc,\Bb)$ is still unknown and can be estimated by replacing
$\boldc$ with the corresponding maximum likelihood estimator. We  denote the
estimated $\Bu_{11} (\boldc,\Bb)$ by $\hat{\Bu}_{11} (\boldc,\Bb)$. Since
$\hat{\boldc}$ is a consistent estimator of $\boldc$, $\hat{\Bu}_{11} (\boldc,\Bb)$
converges in probability to $\Bu_{11} (\boldc,\Bb)$ and therefore a test statistic
for testing the hypothesis \myeq{eq:6d1} is given by
\begin{equation} \label{eq:6d2}
n_{+}(\hat{c}_{2} - 1,...,\hat{c}_{K} - 1)\hat{\Bu}_{11} (\boldc,\Bb) (\hat{c}_{2} -
1,...,\hat{c}_{K} - 1)^{T}
\end{equation}
 whose asymptotic distribution under \myeq{eq:6d1} is a
chi-squared distribution with  $K-1$ degrees of freedom (Rao, 1973, Chapter 6).
 The test statistic \myeq{eq:6d2} can be expressed
in a simple form
$$
n_{+}\{\frac{pr_{1}(1-r_{1})}{2} - pr_{1}\sum_{k=2}^{K}\frac{r_{k}}{\hat{c}_{k}} + \frac{p}{2}
[
\sum_{k=2}^{K}\frac{r_{k}}{\hat{c}_{k}^{2}} - (\sum_{k=2}^{K}
\frac{r_{k}}{\hat{c}_{k}})^{2}] \}$$
$$
 =  \frac{n_{+}p}{2}[\sum_{k=1}^{K}\frac{r_{k}}{\hat{c}_{k}^{2}}
- (\sum_{k=1}^{K}\frac{r_{k}}{\hat{c}_{k}})^{2}],$$
where $\hat{c}_{1} \equiv 1$.

\vspace{.3cm}

\noindent {\bf References} \vspace{.3cm}

\noindent\hangindent .25in
 Anderson, T.W.~(1984).~ {\em An Introduction to Multivariate Statistical Analysis.}  2nd ed.  Wiley, New York.

\noindent\hangindent .25in
Bartlett, M.S.~(1951).~  An inverse matrix adjustment
arising in discriminant  analysis.  {\em Ann. Math. Statist.} {\bf 22}, 107--111.

\noindent\hangindent .25in
 Boente, G., Critchley, F. and  Orellana, L. ~(2007).~ Influence functions of two
families of robust estimators under proportional scatter matrices,  {\em Statistical
Methods and Applications}, {\bf 15}, 295--327

\noindent\hangindent .25in
 Boente, G., Critchley, F. and  Orellana, L. ~(2009).~
Robust tests for the common principal components model, {\em Journal of Statistical
Planning and Inference}, {\bf 139}, 1332--1347

\noindent\hangindent .25in
 Eriksen, P.S.~(1987).~       Proportionality of covariance matrices.
      {\em Ann. Statist.} {\bf 15}, 732--748.

\noindent\hangindent .25in
       Flury, B.N.~(1986).~       Proportionality of k covariance matrices.
       {\em Statist. Prob. Letters} {\bf 4}, 29--33.

\noindent\hangindent .25in
 Guttman, I., D.Y. Kim and I. Olkin~(1985).~
      Statistical inference for constants of proportionality.
      In {\em Multivariate Analysis,}  Vol VI (P.R. Krishnaiah, ed.),
                Elsevier Science Publishers B.V.

\noindent\hangindent .25in
 J\"{o}reskog, K.G.~(1971).~
       Simultaneous factor analysis in several populations.
       {\em Psychometrika} {\bf 36}, 409--426.

\noindent\hangindent .25in
     J\"{o}reskog, K.G. and D. S\"{o}rbom~(1996).~
       {\em LISREL 8: User's reference guide}, Hove and London: Scientific Software
       International.

\noindent\hangindent .25in
 Khatri, C. G.~ (1967). ~
       Some distribution problems connected with the
                 characteristic roots of $S_1 S_2^{-1}$.
       {\em Annals of Mathematical Statistics} {\bf 38}, 944--948.

\noindent\hangindent .25in
       Manly, B. F. and Rayner, J. C. W. (1987).~ The Comparison of Sample Covariance Matrices Using Likelihood
Ratio Tests. {\em Biometrika}, {\bf 74}, 841--847.

\noindent\hangindent .25in
       Owen, A. (1984).~ A neighborhood.based classifier for landsat data. {\em Canadian Journal of Statistics}, {\bf 12}, 191--200.

\noindent\hangindent .25in
       Rao, C.R.~(1973).~
       {\em Linear Statistical Inference and its Applications.}
      2nd ed.
       Wiley, New York.

\noindent\hangindent .25in
       Rao, C. R.~ (1983).~ Likelihood tests for relationships between covariance matrices. In {\em Studies in Econometrics,
Time Series and Multivariate Statistics.}  S. Karlim, T. Ameniya, L. A. Goodman,
editors. pp. 529--543. Academic Press, New York.

\end{document}